*Research Article*
# Some Simple Formulas for Posterior Convergence Rates

**Wenxin Jiang**

*Department of Statistics, Northwestern University, Evanston, IL 60208, USA*

Correspondence should be addressed to Wenxin Jiang; wjiang@northwestern.edu

We derive some simple relations that demonstrate how the posterior convergence rate is related to two driving factors: a "penalized divergence" of the prior, which measures the ability of the prior distribution to propose a nonnegligible set of working models to approximate the true model and a "norm complexity" of the prior, which measures the complexity of the prior support, weighted by the prior probability masses. These formulas are explicit and involve no essential assumptions and are easy to apply. We apply this approach to the case with model averaging and derive some useful oracle inequalities that can optimize the performance adaptively without knowing the true model.

## 1. Introduction

In Jiang [1], there are some general results on the posterior convergence rate, which were very simple and easy to be applied. The current paper is related and developed from the ideas of Jiang [1]. There is no essential new idea behind the proofs. However, the results have been much simplified from the earlier groundbreaking works in this area, such as Ghosal et al. [2], Walker [3], and Ghosal et al. [4] so that the current work can be applied much more easily and displays the intrinsic driving factors behind the convergence rate more directly.

The current paper cannot be used to derive new convergence rates better than what are achievable in the existing literature, except that the current convergence rate is described in a $d_u^2$-divergence for any $u \in (-1, 0)$ (defined later), which is more general than the squared Hellinger distances $d_{-1/2}^2$ (corresponding to $u = -1/2$). A recent work by Norets [5] actually has obtained convergence rates in a stronger Kullback-Leiber divergence $d_0^2$ (corresponding to the limit $u = 0$); however, the rate in $d_0^2$ is suddenly much worse than the rate in $d_{-1/2}^2$ (e.g., about $1/\sqrt{n}$ in parametric cases, instead of about $1/n$). Interestingly, although we allow any $d_u^2$ for $u \in (-1, 0)$, the rate stays at about $1/n$, which is essentially as good as in the Hellinger case at $u = -0.5$ and does not deteriorate to the case of the Kullback-Leibler limit at $u = 0$.

Aside from this technical difference in the divergence measures used in this current paper, the difference from the previous works is essentially esthetical. The key difference from the previous works is that the previous results are presented as *bounds of the posterior probability* (outside a neighborhood of the true density), while the current paper presents almost sure *bounds of the distance* (or *divergences*) from the true density directly. Applications of previous works sometimes need to make a guess on the convergence rate and then check that it simultaneously satisfies several inequality conditions, while the current paper presents explicit formulas that are essentially assumption free.

## 2. Main Results

Let $D$ denote the observed data generated from a probability distribution $P_0$. Consider a prior distribution $\pi$ supported on a set of densities $\mathcal{P}_D$. Define for a subset $A \subset \mathcal{P}_D$ the posterior probability as $\Pi(A) = \int_{f \in A} f(D)\pi(df) / \int_{f \in \mathcal{P}_D} f(D)\pi(df)$. Denote $P_0\Pi(A) = \int \Pi(A) P_0(dD)$. Let $\cup_{j \in N} B_j$ be a countable convex cover of $A$, so that $A \subset \cup_{j \in N} B_j$, $B_j$ is convex, and $N$ is a countable set.

2International Scholarly Research NoticesDefine the divergence $d_t^2(p,q) = t^{-1} \int p((p/q)^t - 1)$ between densities $p$ and $q$, on a suitable dominating measure, for $t > -1$. (The Kullback-Leiber divergence corresponds to the limiting case of $t = 0$; the squared Hellinger distance corresponds to the case $t = -1/2$; the $\chi^2$ divergence corresponds to the case $t = 1$.) Then we have the following result.

**Proposition 1.** *For any $u \in (0,1)$ and any $t \in (0,\infty)$, for any $K \subset \mathcal{P}$, one has*

$$\left(\frac{P_0\Pi(A)}{4}\right)^{1+u/t}$$
$$\leq \sum_{j \in N} \left\{ \frac{\pi(B_j)\left(1 + \sup_{q \in B_j}(-ud_{-u}^2(p_0^D, q^D))\right)^{1/u}}{\pi(K)\left(1 + \sup_{q \in K}(td_t^2(p_0^D, q^D))\right)^{-1/t}} \right\}^u. \quad (1)$$

This result requires no assumption essentially. The relation displayed is explicit. In the result, $p_0^D$, $q^D$ are regarded as probability densities of the entire set of data $D$.

Now consider the case with iid assumption, so that $D = (Y_1, \ldots, Y_n)$, $Y, Y_1, \ldots, Y_n$ being iid (independent and identically distributed), generated from density $p_0$ for a single copy $Y$. Let $\pi$ be a prior distribution supported on a set $\mathcal{P}$ of densities of $Y$. Consider any $A \subset \mathcal{P}$ with any of its countable convex covers $\cup_{j \in \mathcal{N}} B_j$. Using relations such as $p_0^D = \prod_{i=1}^n p_0(Y_i)$ and $1 + x \leq e^x$ for any real $x$, the previous result becomes as follows.

**Proposition 2.** *For iid data with sample size $n$, for any $-1 < -u < 0 < t < \infty$, any $K \subset \mathcal{P}$, one has*

$$\left(\frac{P_0\Pi(A)}{4}\right)^{1+u/t}$$
$$\leq \sum_{j \in \mathcal{N}} \exp\left\{-un\left(\left[\inf_{q \in B_j} d_{-u}^2(p_0, q)\right.\right.\right.$$
$$\left. - n^{-1}\ln\pi(B_j)\right] \quad (2)$$
$$- \left[\sup_{q \in K} d_t^2(p_0, q)\right.$$
$$\left.\left.\left. - n^{-1}\ln\pi(K)\right]\right)\right\}.$$

The only essential assumption here is iid. The relation displayed is explicit.

We will now consider a sequence of densities $p^{(n)}$ for iid data, which are generated from the posterior distributions based on iid data with increasing sample sizes $n$, and study how they converge to the true density $p_0$.

*Condition 1* ("posterior sequence" of random densities for iid data). A "posterior sequence" $p^{(n)}$ (labeled by sample size $n$) of random density functions in $\mathcal{P}$, in a probability space, satisfies, for any subset $A \subset \mathcal{P}$,

$$\Pr\left[p^{(n)} \in A\right] = P_0\Pi(A)$$
$$\equiv \int \left[\frac{\int_{f \in A} \prod_{i=1}^n f(Y_i)\pi(df)}{\int_{f \in \mathcal{P}} \prod_{i=1}^n f(Y_i)\pi(df)}\right] \prod_{i=1}^n p_0(Y_i)\, dY_i. \quad (3)$$

Here $\mathcal{P}$ is a set of density functions, $\pi()$ is the prior distribution of $f \in \mathcal{P}$, and $\Pi()$ is the fraction in the integrand, which can be regarded as the posterior distribution of $f$ based on iid data $D = (Y_1, \ldots, Y_n)$.

At any fixed sample size $n$, this probability law is equivalent to assuming that $p^{(n)}$ is sampled from the posterior $\Pi()$ given data $D = (Y_1, \ldots, Y_n)$, and $D$ is an iid sample of $Y$ with density $p_0$. We will often omit the superscript and write $p^{(n)} = p$.

Suppose $\mathcal{P}$ can be covered by a finite number $\overline{N}$ of the $B_j$'s, each being an $L_1$ ball with radius $n^{-1/u}$. Then the following result can be obtained.

**Proposition 3.** *Consider a "posterior sequence" of densities $p$ for iid data satisfying Condition 1. For any $u \in \{1/2, 1/3, 1/4, \ldots\}$ and any $0 < t < \infty$, with probability 1, for almost all large sample size $n$, one has*

$$d_{-u}^2(p_0, p) \leq d_t^2(p_0, \pi) + n^{-1}\ln\left[\overline{N}^{1/u} n^{2(u^{-1} + t^{-1})}\right], \quad (4)$$

*where $d_t^2(p_0, \pi) \equiv \inf_{K \subset \mathcal{P}}[\sup_{q \in K} d^2(p_0, q) + n^{-1}\ln(1/\pi(K))]$ and $\mathcal{P}$ is the support of the prior $\pi$.*

Here we define a "penalized divergence" $d_t^2(p_0, \pi) \equiv \inf_{K \subset \mathcal{P}}[\sup_{q \in K} d^2(p_0, q) + n^{-1}\ln(1/\pi(K))]$, for some divergence $d$ of a prior $\pi$ from the true density $p_0$.

*Remark 4.* The result can be extended to the continuously valued $u \in (0,1)$. This is because the $d_t$ divergence is monotonically increasing in $t$. For any $d_{-u}$ so that $1/u$ is not an integer, we can use a more stringent divergence $d_{-u'}$ with the $-u'$ being the next larger value from the integer range $\{-1/2, -1/3, -1/4, \ldots\}$ to bound the convergence rate in $d_{-u}$.

*Remark 5.* In this and other works, we notice that we often encounter in the convergence rate results a quantity similar to the "*penalized divergence*" of the form $d_t^2(p_0, \pi) \equiv \inf_K[\sup_{q \in K} d_t^2(p_0, q) + n^{-1}\ln(1/\pi(K))]$, related to a prior $\pi$. This first part $\sup_{q \in K} d_t^2(p_0, q)$ describes the maximal divergence of a set $K$ (proposed by a prior $\pi$) from $p_0$. We can understand this part as the approximation error of the prior $\pi$ when it is used to propose densities to approximate a true density $p_0$. The second part penalizes an unlikely set $K$ with a small prior $\pi(K)$. Combining the two parts, we can perhaps try to interpret $d^2(p_0, \pi)$ as the approximation error (away from $p_0$) by a not-too-unlikely set proposed by a prior $\pi$. This "penalized divergence" is a critically important driving factor for determining the convergence rates in the previous results. It is noted that although this factor corresponds to the approximation ability of $\pi$, it already has a complexity penalty built in it implicitly. This is from the penalty against a small prior; the second part is $n^{-1}\ln(1/\pi(K))$, which is, roughly speaking, about $d/n$, where $d$ is the number of parameters proposed by the prior $\pi$ (e.g., for a uniform prior $\pi$, for a small $d$-dimensional cube $K$ with volume $\delta^d$, we have $\pi(K) \propto \delta^d$).

*Remark 6.* The other factor behind the convergence rate is related to the complexity of the model, which is proportional to $n^{-1}\ln \overline{N}$ where $\overline{N}$ is some number that increases with the number of small convex balls needed to cover the prior support of the model. Typically, this "complexity factor" is roughly about $d/n$, up to some logarithm factors, where $d$ is the dimension of the parameters involved in the prior. It is noted, however, that with model averaging the higher dimensional model can be downweighted by the model prior, so that effectively one can make $d$ to be of order 1 for this complexity factor, so that the convergence rate will be controlled by the first factor ("the penalized divergence") alone.

The convergence rate result in Proposition 3 can be extended to the case of model averaging, when the prior is $\pi(m, df) = \pi_m \pi(df \mid m)$, jointly over a model index $m \in M$ in a set of nonoverlapping models $M$, and density $f \in \mathscr{P}_m$ (the support of prior *given* model $m$) (we assume nonoverlapping models for simplicity, where $\mathscr{P}_m \cap \mathscr{P}_{m'} = \emptyset$ for any two different model indexes $m$ and $m'$. This is only a technical convention for defining the prior supports, which typically does not affect the real applications see, e.g., Section 3.1) and posterior $\Pi(A) \propto \sum_{m\in M} \int_{f\in\mathscr{P}_m} I((m,f) \in A) \prod_{i=1}^n f(Y_i) \pi_m \pi(df \mid m)$ for an event $A$. In this case, let $\overline{N}_m = N((1/n)^{1/u}, L_1, \mathscr{P}_m)$, the $L_1$ balls of radius $(1/n)^{1/u}$ needed to cover the prior-support $\mathscr{P}_m$ under model $m$. Then we have the following, under the iid assumption.

**Proposition 7.** *Consider a "posterior sequence" of densities $p$ for iid data satisfying Condition 1. For any $u \in \{1/2, 1/3, 1/4, \ldots\}$ and any $0 < t < \infty$, with probability 1, for almost all large sample size $n$, one has*

$$d_{-u}^2(p_0, p) \le d_t^2(p_0, \pi) + n^{-1}\ln\left(n^{2(u^{-1}+t^{-1})}\left(\sum_{m\in M}\pi_m^u \overline{N}_m\right)^{1/u}\right), \quad (5)$$

*where*

$$d_t^2(p_0, \pi) \equiv \inf_{K\in\mathscr{P}}\left\{\sup_{q\in K} d_t^2(p_0, q) - n^{-1}\ln \pi(K)\right\}$$

$$\le \inf_{m\in M, K_m\in\mathscr{P}_m}\left\{\sup_{q\in K_m} d_t^2(p_0, q) - n^{-1}\ln[\pi_m \pi(K_m \mid m)]\right\}$$

$$= \inf_{m\in M}\left[d_t^2(p_0, \pi(\cdot \mid m)) - n^{-1}\ln \pi_m\right]. \quad (6)$$

$\mathscr{P}$ and $\mathscr{P}_m$ are supports of the mixing prior $\pi = \sum_m \pi_m \pi(df \mid m)$ and the model-$m$ prior $\pi(df \mid m)$, respectively.

This is an oracle inequality that achieves the best performance of all models $m$ for the bound on the right hand side. Again, the convergence rate is displayed explicitly, and we will try to explain the driving factors of the convergence rate later. This is unlike the previous works where one has to conjecture a rate $\epsilon$ and check that it satisfies many conditions.

So far, we have assumed existence of a finite covering number for the prior support, such as $\overline{N}$ in Proposition 3 or $\overline{N}_m$ in Proposition 7. They determine the "complexity factor" as commented in Remark 6. A deeper analysis of the "complexity factor" is to regard it as an upper bound for a better complexity measure related to the prior $\pi$, developing an idea pioneered by Walker [3].

*Remark 8.* The complexity in Remark 6 is not satisfactory when the prior support $\mathscr{P}$ is unbounded and the covering number $\overline{N}$ is infinity. However, the proofs of the propositions can be easily adapted to show that the covering number $\overline{N}$ can be replaced by $\sum_{j\in\mathcal{N}}\pi(B_j)^u$ from Proposition 2, where we have relaxed $\cup_j B_j$ to be a cover of the entire prior support of $\mathscr{P}$, and we have freedom in choosing the cover $\cup_j B_j$. Therefore, we can define a quantity that is related to the prior itself. Let the $N_u(\pi)$ be the infimum of the $\ell_u$ norm $[\sum_{j\in\mathcal{N}}\pi(B_j)^u]^{1/u}$ over all such covers $\cup_j B_j$ of $\mathscr{P}$, where each $B_j$ is an $L_1$ ball of radius $n^{-1/u}$. We may name it as the "$\ell_u$-norm prior complexity" for covering the prior support. An unbounded prior support may still be coverable by infinitely many $B_j$'s, so that $N_u(\pi)$ is finite, even with an infinite covering number $\overline{N}$. Then we have a better way of formulating a bound corresponding to Proposition 3:

$$d_{-u}^2(p_0, p) \le d_t^2(p_0, \pi) + n^{-1}\ln\left[N_u(\pi) n^{2(u^{-1}+t^{-1})}\right], \quad (7)$$

where $N_u(\pi)$ is the "$\ell_u$-norm complexity" of this prior $\pi$ defined in this remark.

*Remark 9.* We now describe heuristically how to bound the "norm complexity" $N_u(\pi)$ defined in the previous remark in parametric models, where densities $p \in \mathscr{P}$ are parameterized by a $d$ dimensional parameter $\theta$, and a prior on $\theta$ induces a prior $\pi$ on the densities in $\mathscr{P}$. A more rigorous treatment is given in the example of Section 3.2. In typical situations with some smoothness conditions on the densities, we can relate the $L_1$ distance between two densities $p_{\theta_1}$ and $p_{\theta_2}$ by the maximal norm $|\cdot|_\infty$: $\int |p_1 - p_2| \le cd|\theta_1 - \theta_2|_\infty$ for some constant $c > 1$. Then, to cover the parameter space, we can use $\ell_\infty$ ball $A_j$'s in the parameter space with radius $h = (cdn^{1/u})^{-1}$, so that the corresponding densities cover the $L_1$-ball $B_i$ with the required radius $n^{-1/u}$. These sets $A_j$, with small volumes $\text{vol}(A_j) = (2h)^d$, can be used to form a fine partition of the parameter space, so that the norm $\sum_{j\in N}\pi(B_j)^u = \sum_{j\in N}\pi(\theta \in A_j)^u = [\sum_{j\in N}\pi(\theta_j)^u \text{vol}(A_j)]\text{vol}(A_j)^{u-1} \approx [\int \pi(\theta)^u d\theta]((2h)^d)^{u-1}$, where $\pi(\theta_j)$ is the prior density function evaluated at some intermediate point in the set $A_j$. The sum in the square bracket is a Riemann sum over a fine grid, which we will assume to be approximated by an integral under some regularity conditions, even if the domain may be unbounded. Therefore, we have an upper bound of the norm complexity as

$$N_u(\pi) \lesssim \left(\int \pi(\theta)^u d\theta\right)^{1/u} \left(\frac{cdn^{1/u}}{2}\right)^{d(1-u)/u} \tag{8}$$
$$\leq |\pi|_u (cdn)^{d/u^2},$$

for all large enough $n$. Assume that the prior density is $L_u$ integrable in the parameter space, and the norm $|\pi|_u$ scales as $(\text{const})^d$ as in the case of an iid prior $\pi(\theta) = \prod_{j=1}^d \pi_1(\theta_j)$. Then the complexity term in the bound of Remark 8 can be derived as

$$n^{-1} \ln\left[N_u(\pi) n^{2(u^{-1}+t^{-1})}\right] = O\left(n^{-1} d \ln(nd)\right), \tag{9}$$

which increases with the dimension $d$.

*Remark 10.* Similar to Remark 8, we have a better way of formulating a bound corresponding to Proposition 7:

$$d_{-u}^2(p_0, p) \leq d_t^2(p_0, \pi) + n^{-1} \ln\left[N_u(\pi) n^{2(u^{-1}+t^{-1})}\right], \tag{10}$$

where $N_u(\pi)$ is the "$\ell_u$-norm complexity" of this prior $\pi$, which in this case should be the infimum of the $\ell_u$ norm $[\sum_{m \in M} \sum_{j \in \mathcal{N}_m} \{\pi_m \pi(B_j \mid m)\}^u]^{1/u}$ over all such covers $\cup_{m,j} B_{mj}$ of $\mathcal{P}$, where each $B_{mj}$ is an $L_1$ ball of radius $n^{-1/u}$, and under each model $m$, $\cup_{j \in \mathcal{N}_m} B_{mj}$ represents a cover of its prior support $\mathcal{P}_m$ using possibly infinitely many balls. The defining expression of $N_u(\pi)$ can also be related to the norm complexities of all the conditional priors $\pi(\cdot, m)$ given the model choices: $N_u(\pi) = [\sum_{m \in M} \{\pi_m N_u(\pi(\cdot \mid m))\}^u]^{1/u}$. With model averaging using some suitable weights $\pi_m$, this term $N_u(\pi)$ and its effect on the convergence rate no longer diverge with the complexity of the model, in contrast to the conclusion of Remark 9. The convergence rate is then mainly determined by the penalized divergence $d_t^2(p_0, \pi)$. An example below (in its second part) is used to illustrate this.

## 3. A Simple Example for Illustration

This is a simple binary regression example intended for illustration. We will see that model averaging can be used to derive nearly optimal convergence rates that are adaptive to the assumptions on the true model. In the first part, we will illustrate how to bound the penalized divergence with a uniform prior with a bounded support. In the second part, we will illustrate how to bound the norm complexity when the prior has an unbounded support.

*3.1. When the Prior Has a Bounded Support.* Consider a binary regression model $Z \mid X \sim \text{Bin}(1, \mu(X))$, $X \in \text{Unif}(0,1)$, where the true conditional mean function $\mu(X)$ is denoted as $\mu_0(X) \in (\delta, 1-\delta)$ (for some small positive $\delta$), which is bounded away from 0 and 1, for any value of $X$. We consider an $m$-piecewise constant working model for the mean function $\mu$. Suppose the prior is $\pi = \pi_m \pi(\theta_1, \ldots, \theta_m \mid m)$ where $m$ indicates the $m$-piecewise constant model $\mu(x) = \sum_{j=1}^m \theta_j I[x \in [(j-1)/m, j/m)]$. We consider an independent uniform prior $\pi(\theta_1, \ldots, \theta_m \mid m) = \prod_{j=1}^m I(\theta_j \in [0,1])$. (For technically defining the prior supports to be nonoverlapping for different models, one can further require the $\theta_j$'s to be mutually distinct. The resulting prior would be unchanged almost everywhere and would not affect the discussions later.)

We will consider two different setups of the true model.

*Setup 1* (dense true model). In the first setup, the true $\mu_0$ has continuous derivative bounded by $D$. We call this a "dense" setup since we may need a large piecewise constant model (with large $m$ increasing with sample size $n$) to approximate this quite arbitrary true mean function $\mu_0$.

To apply Proposition 7, we will use $-u = -1/2$ and $t = 1$. In the present case, we have $d_1^2(p_1, p_2) = \int_0^1 dx (\mu_1 - \mu_2)^2/(\mu_2(1-\mu_2)), \leq (\delta(1-\delta))^{-1} \int_0^1 dx (\mu_1 - \mu_2)^2$, if $p_1(z, x) = \mu_1^z(1-\mu_1)^{1-z}$ and $p_2(z, x) = \mu_2^y(1-\mu_2)^{1-z}$. We will sometimes use the mean function $\mu$'s to denote the corresponding distances as $d_1^2(\mu_1, \mu_2)$ and so on.

The approximation property of the piecewise constants implies that, for any true $\mu_0(x)$, there exists a $\mu^*$ in the support of $\pi(\cdot \mid m)$ so that $|\mu_0 - \mu^*| < D/m$ everywhere. In fact, we can take $\mu^* = \sum_{j=1}^m \theta_j^* I[x \in [(j-1)/m, j/m)]$ where $\theta_j^* = \mu_0((j-1)/m) \in (\delta, 1-\delta)$. Now let $\mu$ be close to $\mu^*$, in the sense that $\mu = \sum_{j=1}^m \theta_j I[x \in [(j-1)/m, j/m)]$ and $\theta_j \in \theta_j^* \pm \Delta$, $j = 1, \ldots, m$, for some small $\Delta \in (0, \delta)$. Take $K_m$ to be the set of all such densities of $\text{Bin}(1, \mu)$. Then $\mu \in (\delta - \Delta, 1 - \delta + \Delta)$ everywhere, since $\theta_j \in \theta_j^* \pm \Delta$ and $\theta_j^* \in (\delta, 1-\delta)$ for all $j$.

Let $p_0$ and $p$ be the densities corresponding to $\mu_0$ and $\mu$, respectively. Then $\sup_{p \in K_m} d_1^2(p_0, p) = \sup_{p \in K_m} \int_0^1 dx(\mu_0 - \mu)^2/(\mu(1-\mu)) \leq \sup_{p \in K_m} \int_0^1 dx(|\mu_0 - \mu^*| + |\mu - \mu^*|)^2/\{\min_x(\mu(1-\mu))\} \leq (D/m + \Delta)^2/((\delta - \Delta)(1 - \delta + \Delta))$ due to the triangle inequality. The prior probability over $(m, K_m)$ is $\pi_m(2\Delta)^m$. Therefore, the "penalized divergence"

$$d_1^2(p_0, \pi) \leq \inf_{m, \Delta} \left\{\left((\delta - \Delta)(1 - \delta + \Delta)^{-1}\right)\left(\frac{D}{m} + \Delta\right)^2 + \left(\frac{m}{n}\right) \ln(2\Delta)^{-1} + n^{-1} \ln \pi_m^{-1}\right\}. \tag{11}$$

We will take $\pi_m \propto e^{-K(m-1)\ln n}$ for some large enough constant $K > 0$ and apply Proposition 7. (It can be shown that this will make the "complexity" term $n^{-1} \ln(n^6(\sum_{m=1}^\infty \sqrt{\pi_m} \overline{N}_m)^2)$ negligible compared to $d^2(p_0, \pi)$ (by showing $\overline{N}_m \leq n^{2m}$, we omit the tedious details here).)

We can take $m \sim n^{1/3}$ and $\Delta = 1/n$ for an upper bound of the $\inf_{m,\Delta}$. Therefore, the "penalized divergence" $d_1^2(p_0, \pi)$ and the resulting convergence rate $d_{-1/2}^2(p_0, p)$ are both of order $O(n^{-2/3} \ln n)$, which is within a $\ln n$ factor to the minimax optimal result. It is noted that the model averaging automatically achieves this near optimal rate.

*Setup 2* (sparse true model). Consider a second setup, where we assume that the true model is a $m_0$-piecewise constant, where we do not know the value $m_0$. We call this a sparse case since we only need an $m$-piecewise constant model $\mu^*$ to approximate the true mean function $\mu_0$ perfectly, where $m = m_0$ can be much smaller than the choice of $m \sim n^{1/3}$ in Setup 1.

Then modifying the above reasoning, we can bound the infimum in the "penalized divergence" $d_1^2(p_0, \pi)$ by taking $m = m_0$ and $\Delta = 1/n$ and obtain (using $|\mu^* - \mu_0| = 0$) the following:

$$d_1^2(p_0, \pi) \leq \left(\left(\delta - \frac{1}{n}\right)\left(1 - \delta + \frac{1}{n}\right)\right)^{-1}\left(0 + \frac{1}{n}\right)^2 \\ + O\left(\frac{m_0 \ln n}{n}\right) = O\left(\frac{m_0 \ln n}{n}\right). \quad (12)$$

This and also the resulting posterior convergence rate $d_{-1/2}^2(p_0, p)$ are, therefore, both close to the parametric rate $O(m_0/n)$, as if we knew $m_0$ beforehand.

In summary, the prior $\pi$ is *adaptive* in the sense that in either the dense or the sparse case, the resulting posterior distribution works nearly optimally, even if we do not really know whether the true model is dense or sparse.

*3.2. When the Prior Has an Unbounded Support.* In the example above, we have considered a uniform prior for $\pi(\cdot \mid m)$ with a bounded support. In this subsection, we will consider a parametrization in the log-odds scale, with an unbounded prior support, for illustrating how to calculate the norm complexity described in Remarks 8, 9, and 10. The model is still $Z \mid X \sim \text{Bin}(1, \mu)$, $X \in \text{Unif}(0, 1)$, where the true $\mu$, denoted as $\mu_0$, is bounded away from 0 and 1. We consider an $m$-piecewise constant working model for the mean function $\mu = \sum_{j=1}^m e^{\theta_j}/(1 + e^{\theta_j})I[x \in [(j-1)/m, j/m]]$. Suppose the prior is iid for each "log-odds" parameter $\theta_j$ supported on $\Re$; then

$$\pi = \pi_m \pi(\theta_1, \ldots, \theta_m \mid m) = \pi_m \prod_{j=1}^m \pi(\theta_j \mid m). \quad (13)$$

Consider two densities of $(Z, X)$: $p_{1,2} = \mu_{1,2}^z (1 - \mu_{1,2})^{1-z}$, with parameters $\theta_{1j}$'s and $\theta_{2j}$'s, respectively. Then one can easily derive a following relationship for the $L_1$ distance: $\int |p_1 - p_2| = 2\int |\mu_1 - \mu_2| \leq (1/2)\sup_{j=1}^m |\theta_{1j} - \theta_{2j}|$. For covering all the densities in this working model with $L_1$ balls $B_{mj}$ with radius $n^{-1/u}$ ($u \in (0, 1)$), we can use the densities with parameters in $\ell_\infty$ balls $[\theta_s \in [j_s h, (j_s + 1)h] : s = 1, \ldots, m]$, with radius $h/2 = 2n^{-1/u}$.

Then $N_u(\pi(\cdot \mid m)) \leq [\sum_{j_1=-\infty}^\infty \cdots \sum_{j_m=-\infty}^\infty \pi(\theta_1 \in [j_1 h, (j_1+1)h], \ldots, \theta_m \in [j_m h, (j_m+1)h] \mid m)^u]^{1/u} = [\prod_{s=1}^m \sum_{j_s=-\infty}^\infty \pi(\theta_s \in [j_s h, (j_s + 1)h] \mid m)^u]^{1/u} =$
$[\sum_{j=-\infty}^\infty \pi(\theta_1 \in [jh, (j+1)h] \mid m)^u]^{m/u}$, since the priors of $\theta_1, \ldots, \theta_m \mid m$ are iid.

Now assume that the prior density for each $\theta_j$ (given model $m$) is continuous, symmetric, and decreasing from the origin $\pi(\theta_j \mid m) = f_m(|\theta_j|)$ for some decreasing functions $f_m \in L_u(0, \infty)$ (which can be satisfied by, e.g., independent $N(0, 1)$ priors or double-exponential densities).

Then $\sum_{j=-\infty}^\infty \pi(\theta_1 \in [jh, (j+1)h] \mid m)^u = 2\sum_{j=0}^\infty \pi(\theta_1 \in [jh, (j+1)h] \mid m)^u \leq 2\sum_{j=0}^\infty f_m(jh)^u h^{u-1} \leq [2hf_m(0)^u + \int_{-\infty}^\infty f_m(|w|)^u dw]h^{u-1}$, where we used the decreasingness of $f_m$ in the last two steps. The integration exists since $f_m \in L_u(0, \infty)$.

Then we have $N_u(\pi(\cdot \mid m)) \leq O(h^{m(u-1)/u}) = O(n^{m(1-u)/u^2})$, which is finite despite the unbounded priors support. Then according to Remark 10, $N_u(\pi) = [\sum_m (\pi_m N_u(\pi(\cdot \mid m)))^u]^{1/u}$ will be $O(n^{(1-u)/u^2})$ if $\pi_m \propto e^{-K(m-1)\ln n}$ for some large enough constant $K > 0$. So the norm complexity term in Remark 10 is of order $O(\ln n/n)$, which, when compared with the last formula in Remark 9, behaves as if the dimension has become reduced to order $O(1)$ by model averaging. Therefore, the norm complexity term does not affect the convergence rate significantly due to model averaging, and the convergence rate is mainly determined by the penalized divergence $d_1^2(p_0, \pi)$. The bounding of the penalized divergence is similar to the example discussed in the previous subsection and we omit the details. The resulting convergence rates are essentially the same as when the uniform priors (with bounded supports) are used, despite the fact that we now allow priors with unbounded supports (such as normal priors in the parametrization of log-odds).

## 4. Proofs

*Proof of Proposition 1.* $P_0 \Pi(A) \leq P_0 \phi + P_0 \Pi(A)(1 - \phi)$ for any "test" $\phi$ as a function of $D$ valued in $[0, 1]$. Consider

$$P_0 \Pi(A)(1 - \phi) \\ = P_0 \left[ \frac{\int_{f \in A} (f(D)/p_0(D))(1 - \phi)\pi(df)}{\int (f(D)/p_0(D))\pi(df)} \right] \\ \leq \frac{P_0 \int_{f \in A} (f(D)/p_0(D))(1 - \phi)\pi(df)}{a} \\ + P_0 I\left[\int \left(\frac{f(D)}{p_0(D)}\right)\pi(df) \leq a\right] \quad \text{for any } a > 0, \quad (14)$$

where $P_0$ represents the expectation under the true density $p_0$.

$P_0 \int_{f \in A} (f(D)/p_0(D))(1 - \phi)\pi(df) = \int_{f \in A} P_f(D)(1 - \phi)\pi(df)$ due to Fubini's Theorem, where $P_f$ represents the expectation under density $f$.

Using Markov's theorem, for any $t > 0$, we have

$$P_0 I\left[\int\left(\frac{f(D)}{p_0(D)}\right)\pi(df) \le a\right]$$

$$\le a^t P_0\left(\frac{\int f(D)\pi(df)}{p_0(D)}\right)^{-t}$$

$$\le a^t P_0\left(\frac{\int_{f\in K} f(D)\pi(df)}{p_0(D)}\right)^{-t},$$

for any $K$ in the support of $\pi$.

$$= a^t \pi(K)^{-t} P_0\left(\frac{\int_{f\in K}(f(D)/p_0(D))\pi(df)}{\pi(K)}\right)^{-t}$$

$$\le a^t \pi(K)^{-t} P_0\left(\frac{\int_{f\in K}(f(D)/p_0(D))^{-t}\pi(df)}{\pi(K)}\right) \quad (15)$$

due to Jensen's inequality;

$$= a^t \pi(K)^{-t} \left(\frac{\int_{f\in K} P_0(f(D)/p_0(D))^{-t}\pi(df)}{\pi(K)}\right)$$

due to Fubini's Theorem;

$$\le a^t \pi(K)^{-t} \sup_{f\in K} P_0\left(\frac{f(D)}{p_0(D)}\right)^{-t}.$$

All these combine to $(*)$

$$P_0\Pi(A) \le P_0\phi + \frac{\int_{f\in A} P_f(D)(1-\phi)\pi(df)}{a}$$
$$+ a^t \pi(K)^{-t} \sup_{f\in K} P_0\left(\frac{f(D)}{p_0(D)}\right)^{-t}. \quad (16)$$

Now apply a result that is a straightforward extension of Ghosal et al. ([4], Lemma 6.1). For any convex set $B_j$, there exist $\phi_j$ such that for any $\alpha, \beta > 0$, $u \in (0,1)$, $\sup_{f\in B_j}(\alpha P_0\phi_j + \beta P_f(1-\phi_j)) = \sup_{f\in B_j}\int\min\{(\alpha p_0),(\beta f)\} \le \sup_{f\in B_j}\int\{(\alpha p_0)^{1-u}(\beta f)^u\} = \alpha^{1-u}\beta^u \sup_{f\in B_j} P_0(p_0/f)^{-u}$.

Therefore, we can find $\phi_j$ so that

$$P_0\phi_j \le \alpha^{-u}\beta^u \sup_{f\in B_j} P_0\left(\frac{p_0}{f}\right)^{-u} \equiv \lambda_j^{-u}\sup_{f\in B_j} P_0\left(\frac{p_0}{f}\right)^{-u},$$

$$\sup_{f\in B_j} P_f(1-\phi_j) \le \alpha^{1-u}\beta^{u-1}\sup_{f\in B_j} P_0\left(\frac{p_0}{f}\right)^{-u} \quad (17)$$

$$= \lambda_j^{1-u}\sup_{f\in B_j} P_0\left(\frac{p_0}{f}\right)^{-u}.$$

Given any $\lambda_j > 0$, we can choose $(\alpha, \beta)$ so that $\alpha/\beta = \lambda_j$ in the above statement.

If we take $\cup_{j\in N} B_j$ to be a convex cover of $A$ and define a combined test $\phi = \max_{j\in N}\phi_j$ and plug it into $(*)$, we have

$$P_0\Pi(A)$$

$$\le \sum_{j\in N} P_0\phi_j + \frac{\sum_{j\in N}\int_{f\in B_j} P_f(D)(1-\phi_j)\pi(df)}{a}$$

$$+ a^t \pi(K)^{-t}\sup_{f\in K} P_0\left(\frac{f(D)}{p_0(D)}\right)^{-t}$$

$$\le \sum_{j\in N} \lambda_j^{-u} \sup_{f\in B_j} P_0\left(\frac{p_0}{f}\right)^{-u} \quad (18)$$

$$+ \frac{\sum_{j\in N}\pi(B_j)\lambda_j^{1-u}\sup_{f\in B_j} P_0(p_0/f)^{-u}}{a}$$

$$+ a^t\pi(K)^{-t}\sup_{f\in K} P_0\left(\frac{f(D)}{p_0(D)}\right)^{-t}.$$

Therefore, we have

$$P_0\Pi(A) \le \sum_{j\in N}\left(\lambda_j^{-u} + \frac{\pi(B_j)\lambda_j^{1-u}}{a}\right)\sup_{f\in B_j} P_0\left(\frac{p_0}{f}\right)^{-u} \quad (19)$$

$$+ a^t\pi(K)^{-t}\sup_{f\in K} P_0\left(\frac{f(D)}{p_0(D)}\right)^{-t}.$$

Take $\lambda_j = (1-u)^{-1} u a \pi(B_j)^{-1}$ and

$$a = \left[\frac{\sum_{j\in N}\pi(B_j)^u \sup_{f\in B_j} P_0(p_0/f)^{-u}}{\{\pi(K)^{-t}\sup_{f\in K} P_0(p_0/f)^t\}}\right]^{1/(u+t)} \quad (20)$$

$$\times \left\{[(1-u)^{1-u}u^u]^{-1}ut^{-1}\right\}^{1/(u+t)}.$$

We get

$$P_0\Pi(A)$$

$$\le c(u,t)\left(\sum_{j\in N}\left\{\frac{\pi(B_j)\left[\sup_{f\in B_j} P_0(p_0/f)^{-u}\right]^{1/u}}{\pi(K)\left[\sup_{f\in K} P_0(p_0/f)^t\right]^{-1/t}}\right\}^u\right)^{t/(u+t)}, \quad (21)$$

where $c(u,t) = (t/(t+u))^{-t/(t+u)}(u/(t+u))^{-u/(t+u)}(u^u(1-u)^{1-u})^{-t/(u+t)} \le 4$. Then notice that $P_0(p_0/f)^t = 1 + td_t^2(p_0, f)$; we obtain the proof of Proposition 1. (For notational convenience, the densities which appeared in the proof here, such as $p_0$ and $f$, are the densities for the entire data set $p_0^D$ and $f^D$, resp., and we do not assume iid.) □

*Proof of Proposition 2.* Use the fact that, for any $t \in (-1, 0) \cup (0, \infty)$,

$$\begin{aligned}
& 1 + t d_t^2 \left(p_0^D, f^D\right) \\
&= P_0 \left(\frac{p_0^D}{f^D}\right)^t = \int p_0^D \left(\frac{p_0^D}{f^D}\right)^t \\
&= \left(\int p_0 \left(\frac{p_0}{f}\right)^t\right)^n, \text{ due to the iid assumption.} \\
&= \left(1 + t d_t^2 (p_0, f)\right)^n \leq e^{t n d_t^2 (p_0, f)}.
\end{aligned} \quad (22)$$

This leads to the proof by applying Proposition 1. □

*Proof of Proposition 3.* This is a special case of Proposition 7, where $\pi_m$ focuses on only one model. □

*Proof of Proposition 7.* Repeat the proofs of Propositions 1 and 2 for the case with model averaging, with the support of prior being $\mathscr{P} = \cup_{m \in M} \mathscr{P}_m$, where $\mathscr{P}_m$ is the support of $\pi(\cdot \mid m)$.

Suppose the convex cover of $A \cap \mathscr{P}$ is doubly indexed as $\cup_{m \in M} \cup_{j \in \mathscr{N}_m} B_{mj}$, where $\mathscr{N}_m$ has cardinality at most $\overline{N}_m$ and $\cup_{j \in \mathscr{N}_m} B_{mj}$ is a convex cover of the support $A \cap \mathscr{P}_m$. Then the result in Proposition 2 holds with

$$\begin{aligned}
& \left(\frac{P_0 \Pi(A)}{4}\right)^{1+u/t} \\
& \leq \sum_{m \in M, j \in \mathscr{N}_m} \exp\left\{-u\left(\left[n \inf_{q \in B_{mj}} d_u^2(p_0, q)\right.\right.\right. \\
&\qquad\qquad\qquad\qquad - \ln\left[\pi_m \pi(B_{mj} \mid m)\right]\right] \\
&\qquad\qquad\qquad\qquad \left.\left.- \left[n \sup_{q \in K} d_t^2(p_0, q) - \ln \pi(K)\right]\right)\right\}.
\end{aligned} \quad (23)$$

In Proposition 2, let $A = [f : d_{-u}^2(p_0, f) > \epsilon]$.

Suppose all the convex sets $B_{mj}$ are such that $\inf_{j,m} \inf_{f \in B_{mj}} d_{-u}^2(p_0, f) > \delta$. Then we have (†)

$$\begin{aligned}
& \left(\frac{P_0 \Pi(A)}{4}\right)^{1+u/t} \\
& \leq \sum_{m \in M} \pi_m^u \overline{N}_m \\
& \quad \times \exp\left\{-un\left(\delta - \left[\sup_{q \in K} d_t^2(p_0, q) - n^{-1} \ln \pi(K)\right]\right)\right\},
\end{aligned} \quad (24)$$

where $\overline{N}_m$ is an upper bound for the number of convex sets $B_{mj}$ needed to cover $A \cap \mathscr{P}_m$.

Now we try to define the convex sets $B_{mj}$ in more detail. They are used to cover $A$, so without generality, each $B_{mj}$ contains a point in $A$, say $p_1$, which is not close to $p_0$ since $d_{-u}^2(p_0, p_1) > \epsilon$. If $B_{mj}$ is small so that any two points are close together, then any point $p_2$ in $B_{mj}$ (which may fall outside $A$) can be made to be also not close to $p_0$, so that $d_{-u}^2(p_0, p_2) > \delta$ for some $\delta > 0$ related to $\epsilon$. This would be easy to establish by a triangular inequality, were it not for the difficulty that the divergence $d_{-u}$ is not a true distance for $u \neq 1/2$. So we would not be able to say, for example, that $B_{mj}$ should be a small $d_{-u}$-ball.

To resolve this difficulty, we derive the following inequalities:

$$\begin{aligned}
& u \left|d_{-u}^2(p_0, p_1) - d_{-u}^2(p_0, p_2)\right| \\
& \leq \int p_0^{1-u} |p_1^u - p_2^u| \leq \left\{\int |p_1^u - p_2^u|^{1/u}\right\}^u
\end{aligned} \quad (25)$$

due to Hölder's inequality, and for $u = 1/k$ for any $k = 2, 3, 4, \ldots$, we have $|p_1^u - p_2^u|^{1/u} = |p_1^{1/k} - p_2^{1/k}||p_1^{1/k} - p_2^{1/k}|^{k-1} \leq |p_1^{1/k} - p_2^{1/k}| \sum_{i=0}^{k-1} |p_1^{i/k} p_2^{(k-1-i)/k}| 2^{k-1} = |(p_1^{1/k})^k - (p_2^{1/k})^k| 2^{k-1}$, so that

$$\left\{\int |p_1^u - p_2^u|^{1/u}\right\}^u \leq 2 \left(\int |p_1 - p_2|\right)^u. \quad (26)$$

Therefore, if $B_{mj}$ is an $L_1$ ball with small enough radius $\lambda$, then any point $p_2$ in $B_{mj}$ has a not too small divergence $d_{-u}^2(p_0, p_2) \geq \epsilon - |d_{-u}^2(p_0, p_1) - d_{-u}^2(p_0, p_2)| \geq \epsilon - u^{-1} 2 (2\lambda)^u \geq \epsilon - u^{-1} 4 \lambda^u = \delta$. We can take $\lambda = [(u/4)(\epsilon - \delta)]^{1/u}$. The $L_1$ balls are also convex as required. We will take $\epsilon - \delta = 4/(nu)$. Then the radius $\lambda = (1/n)^{1/u}$. Therefore, we conclude that we will take $B_{mj}$'s to be small $L_1$ balls with radius $n^{-1/u}$.

Now we try to find $\overline{N}_m$, an upper bound of the number of such small balls needed to cover $A \cup \mathscr{P}_m$. We can use this upper bound $\overline{N}_m \equiv N((1/n)^{1/u}, L_1, \mathscr{P}_m)$, the $L_1$ balls of radius $(1/n)^{1/u}$ needed to cover the entire prior-support $\mathscr{P}_m$ under model $m$.

Then (†) implies, for any $\epsilon > 0$, (‡),

$$\begin{aligned}
& \left(\frac{P_0 \Pi\left(d_{-u}^2(p_0, p) > \epsilon\right)}{4}\right)^{1+u/t} \\
& \leq \exp\left\{-un\left(\epsilon - \frac{\ln\left(e^4 \sum_{m \in M} \pi_m^u \overline{N}_m\right)}{(nu)} - d(p_0, \pi)\right)\right\},
\end{aligned} \quad (27)$$

where $d_t^2(p_0 \pi) \equiv \inf_{K \in \mathscr{P}} [\sup_{q \in K} d_t^2(p_0, q) - n^{-1} \ln \pi(K)]$ and $\mathscr{P}$ is the support of the prior.

Let the probability in (‡) be bounded by $4/(e^{-2/(1+u/t)} n)^2$ under a choice of $\epsilon = \epsilon_n$ for the right-hand side of (‡). Then the event $[d_{-u}^2(p_0, p) \leq \epsilon_n]$ will happen for all large enough $n$, almost surely, due to the Borel-Cantelli lemma.

Then

$$(e^{-2/(1+u/t)}n)^{-2(1+u/t)}$$
$$= \exp\left\{-un\left(\epsilon_n - \frac{\ln\left(e^4 \sum_{m\in M} \pi_m^u \overline{N}_m\right)}{(nu)} - d(p_0, \pi)\right)\right\}. \quad (28)$$

This leads to

$$\epsilon_n = \frac{\left[2(1+u/t)\ln n + \ln\left(\sum_{m\in M} \pi_m^u \overline{N}_m\right)\right]}{(nu)} + d_t^2(p_0, \pi). \quad (29)$$

The quantity $d_t^2(p_0, \pi) = \inf_{K\in\mathscr{P}}[\sup_{q\in K} d_t^2(p_0, q) - n^{-1}\ln\pi(K)]$, under a mixture prior $\pi(K) = \sum_{s\in M} \pi_s \pi(K \mid s)$, can be bounded by taking $K = K_m \subset \mathscr{P}_m$, for any $m \in M$, as $\inf_K[\sup_{q\in K} d_t^2(p_0, q) - n^{-1}\ln\pi(K)] \leq \inf_m \inf_{K_m}\{\sup_{q\in K_m} d_t^2(p_0, q) - n^{-1}\ln[\pi_m \pi(K_m \mid m)]\} \equiv \inf_m[d_t^2(p_0, \pi(\cdot \mid m)) - n^{-1}\ln\pi_m]$. □

## Conflict of Interests

The author declares that there is no conflict of interests regarding the publication of this paper.

## Acknowledgments

This work was based on Technical Report 14-01, Department of Statistics, Northwestern University. The author thanks the reviewers for their helpful comments. The author also thanks Qilu Securities Institute for Financial Studies, Shandong University, for the hospitality during his visit, when a revision of this work was done.

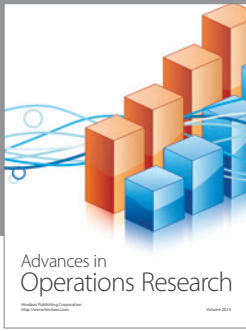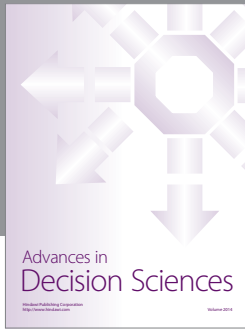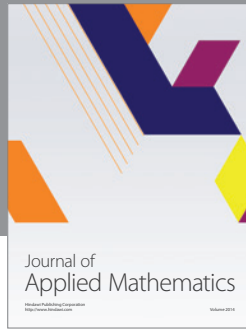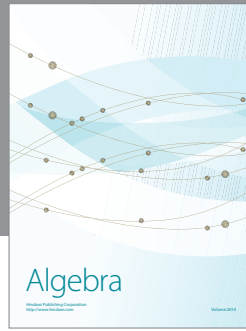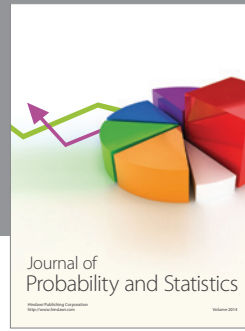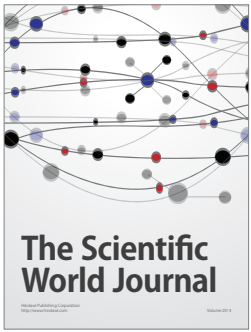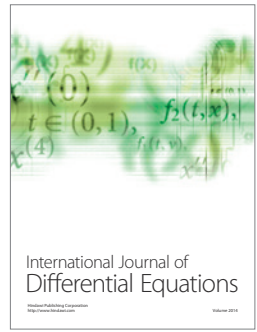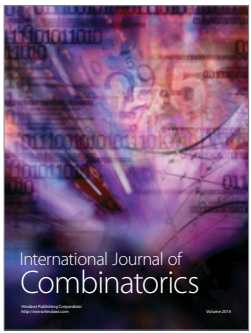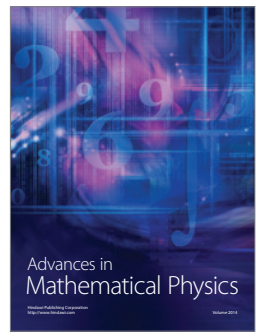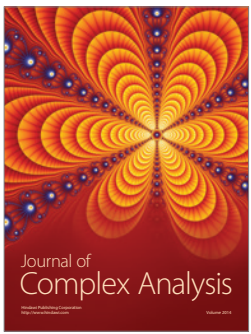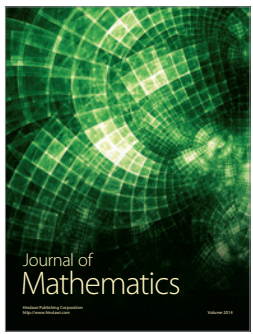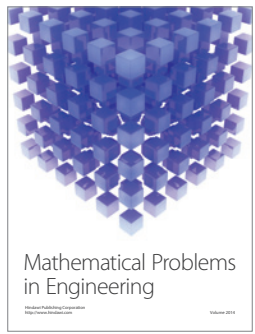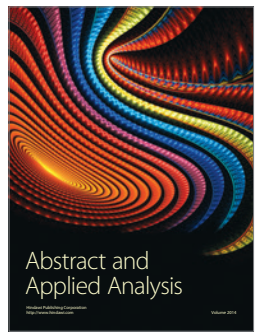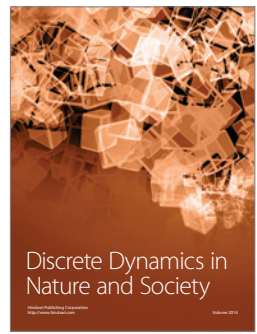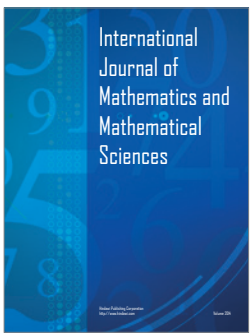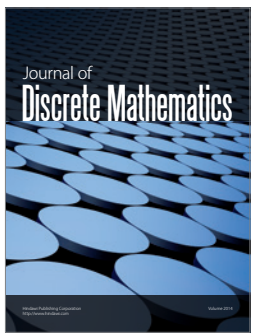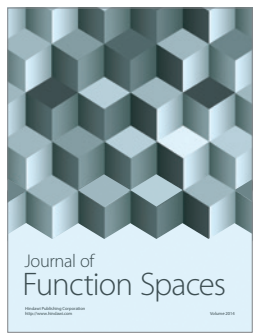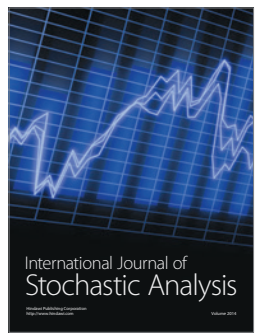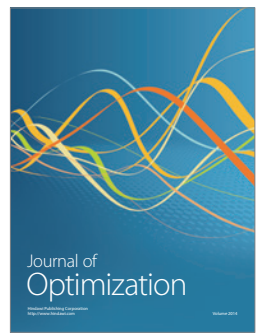

Submit your manuscripts at
http://www.hindawi.com